\documentclass[12pt]{article}
\usepackage{amsthm}
\usepackage{amsmath}
\usepackage{amsfonts}
\usepackage{amssymb}
\usepackage{graphicx}
\usepackage{color}
\usepackage[english]{babel}
\usepackage{soul}
\usepackage[utf8]{inputenc}

\newcommand\blfootnote[1]{%
	\begingroup
	\renewcommand\thefootnote{}\footnote{#1}%
	\addtocounter{footnote}{-1}%
	\endgroup
}

\theoremstyle{plain}   
\newtheorem{theorem}{Theorem}[section]
\newtheorem{proposition}{Proposition}[section]
\newtheorem{corollary}{Corollary}[section]
\newtheorem{lemma}{Lemma}[section]

\DeclareMathOperator{\spec}{sp\,}

\def\R{\ns R}
\def\Z{\ns Z}

\def\vec0{\mbox{\boldmath $0$}}

\def\A{\mbox{\boldmath $A$}}

\def\I{\mbox{\boldmath $I$}}
\def\J{\mbox{\boldmath $J$}}

\def\PP{\mbox{\boldmath $P$}}

\def\I{\mbox{\boldmath $I$}}
\def\J{\mbox{\boldmath $J$}}

\def\R{\mbox{\boldmath $R$}}
\def\Z{\mbox{\boldmath $Z$}}
\def\1{\mbox{\boldmath $1$}}

\newcommand{\tr}[1]{\mathrm{tr\,}#1}
\newcommand{\dist}{\mathrm{dist}}

\title{Almost Moore and the largest mixed \\ graphs of diameters two and three}

\author{
	C. Dalf\'o$^a$, M. A. Fiol$^b$, N. L\'opez$^c$\\\
	\\
	{\small $^a$Dept. de Matem\`atica, Universitat de Lleida, Igualada (Barcelona), Catalonia}\\
	{\small {\tt cristina.dalfo@udl.cat}}\\
	{\small $^{b}$Dept. de Matem\`atiques, Universitat Polit\`ecnica de Catalunya, Barcelona, Catalonia} \\
	{\small Barcelona Graduate School of Mathematics} \\
           {\small Institut de Matem\`atiques de la UPC-BarcelonaTech (IMTech)}\\
           {\small {\tt miguel.angel.fiol@upc.edu}}\\
	{\small $^c$ Dept. de Matem\`atica, Universitat de Lleida, Lleida, Spain}\\
{\small {\tt nacho.lopez@udl.cat}}
\thanks{The research of C. Dalf\'o and M. A. Fiol has been partially supported by
	the grant 2017SGR1087 from AGAUR and the
	grant PGC2018-095471-B-I00 from MCIN/AEI/10.13039/
	501100011033 and “ERDF A way of making Europe", by the European Union. The research of C. Dalf\'o and N. L\'opez has been partially supported by the grant 2017SGR1158 from AGAUR and 
	the grant PID2020-115442RB-I00 from MCIN/AEI/10.13039/501100011033.}}

\date{}
\begin{document}

\maketitle

	\blfootnote{
	\begin{minipage}[l]{0.3\textwidth} \includegraphics[trim=10cm 6cm 10cm 5cm,clip,scale=0.15]{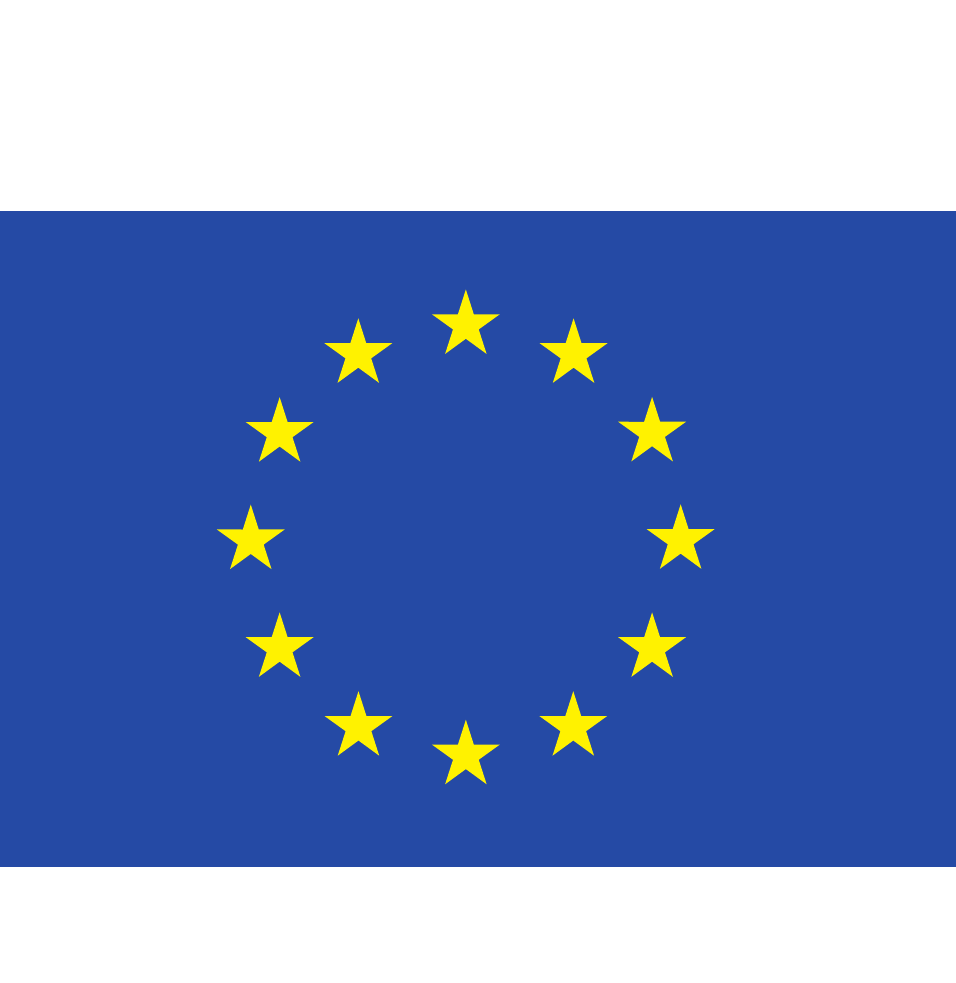} \end{minipage}  \hspace{-2cm} \begin{minipage}[l][1cm]{0.79\textwidth}
		The research of C. Dalf\'o has received funding from the European Union's Horizon 2020 research and innovation program under the Marie Sk\l{}odowska-Curie grant agreement No 734922.
	\end{minipage}
}

\begin{abstract}
{\em Almost Moore mixed graphs\/} appear in the context of the {\em
degree/dia\-meter problem} as a class of extremal mixed graphs,  in the sense that
their order is one unit less than the Moore bound for such graphs. The problem of their existence has been considered just for diameter $2$. In this paper, 
we give a complete characterization of these extremal mixed graphs for diameters 
$2$ and $3$.
We also derive some optimal constructions for other diameters.
\end{abstract}


%
\section{Introduction}
The relationship between vertices or nodes in interconnection networks can be undirected or directed depending on whether the communication between nodes is two-way or only one-way. {\em Mixed graphs} arise in this case and in many other practical situations, where both kinds of connections are needed. Urban street networks are perhaps the most popular ones. Therefore, a {\em mixed graph} $G$ may contain (undirected) {\em edges} as well as directed edges (also known as {\em arcs}). Mixed graphs whose vertices represent the processing elements and whose edges represent their links have been studied before (see Bosák \cite{b79}, Dobravec and Robi\v{c} \cite{dr09}, and Nguyen, Miller, and Gimbert \cite{nmg07}). It is, therefore, natural to consider network topologies based on mixed graphs, and investigate the corresponding degree/diameter problem.

\begin{itemize}
\item {\it Degree/diameter problem for mixed graphs}: Given three natural numbers $r$, $z$,
and $k$, find the largest possible number of vertices $N(r,z,k)$ in a
mixed graph with maximum undirected degree $r$, maximum directed out-degree $z$, and diameter $k$.
\end{itemize}

A natural upper bound for $N(r,z,k)$, known as a Moore(-like) bound, is obtained by counting the number of vertices of a Moore tree $MT(u)$ rooted at a given vertex $u$, with depth equal to the diameter $k$, and assuming that for any vertex $v$, there exists a unique shortest path of length at most $k$ 
from $u$ to $v$. The exact value for this number, which is denoted by $M(r,z,k)$, was given by Buset, El Amiri, Erskine, Miller, and P\'erez-Ros\'es \cite{baemp15} (see also Dalf\'o, Fiol, and L\'opez \cite{dfl18} for an alternative computation), and it turns out to be the following:

\begin{equation}
\label{eq:moorebound4}
M(r,z,k)=A\frac{u_1^{k+1}-1}{u_1-1}+B\frac{u_2^{k+1}-1}{u_2-1},
\end{equation}
where
\begin{align*}
A   &=\displaystyle{\frac{\sqrt{v}-(z+r+1)}{2\sqrt{v}}}, \qquad
B   =\displaystyle{\frac{\sqrt{v}+(z+r+1)}{2\sqrt{v}}},\\
v   &=(z+r)^2+2(z-r)+1, \\
u_1 &=\displaystyle{\frac{z+r-1-\sqrt{v}}{2}}, \qquad
u_2 =\displaystyle{\frac{z+r-1+\sqrt{v}}{2}}.
\end{align*}

This bound applies when $G$ is totally regular with degrees $(r,z)$. In this context, we deal with mixed graphs containing at least one edge and one arc. Mixed graphs of diameter $k$, 
maximum undirected degree $r$ at least 1, maximum out-degree $z$ at least 1, and order given by \eqref{eq:moorebound4} are called {\em mixed Moore graphs}. In the case of diameter 2, such extremal mixed graphs are totally regular of degree $d=r+z$, and they have the property that for any ordered pair $(u,v)$ of vertices, there is a unique walk of length at most $2$ between them. Although some such Moore mixed graphs of diameter two are known to exist, and they are unique (see Nguyen, Miller, and Gimbert \cite{nmg07}), the general problem remains unsettled. Again for diameter 2, Bos\'ak \cite{b79} gave a necessary condition for the existence of a Moore mixed graph, but recently it was proved that there is no Moore mixed graph for the $(r,z)$ pairs $(3,3)$, $(3,4)$, and $(7,2)$ satisfying such necessary condition (see L\'opez, Miret, and Fern\'andez \cite{lmf15}). In general, there are infinitely many pairs $(r,z)$ satisfying Bos\'ak necessary condition for which the existence of a Moore mixed  graph is not known yet. For diameter $k \geq 3$, it was proved that Moore mixed  graphs do not exist, see Nguyen, Miller, and Gimbert \cite{nmg07}.


Because Moore mixed  graphs are quite rare, another line of research focuses on the existence of mixed graphs with prescribed degree and diameter and order just one unit less than the Moore bound. These mixed graphs are known as {\em almost Moore mixed graphs}. They have been extensively studied for the undirected case (Erd\H{o}s, Fajtlowitcz, and Hoffman \cite{efh80}) and for the directed case (Gimbert \cite{g01}). Every almost Moore mixed graph $G$ of diameter $k$ has the property that, for each vertex $v\in V(G)$, there exists only one vertex, denoted by $\sigma(v)$ and called the {\em repeat\/} of $v$, such that there are exactly two walks of length at most $k$ from $v$ to $\sigma(v)$. If $\sigma(v)=v$, then $v$ is called a {\em selfrepeat} vertex. 

\section{Almost Moore mixed graphs of diameter $2$}

\begin{table}[t]
	\begin{center}
		\begin{tabular}{|c|c|c|l|l|c|}
			\hline
			$r$ & $c_1$ & $ c_2$ & $z$ & $n$ & Existence \\ \hline
			4   & - & 3 & $1,4,7,10,\dots$ & $26,68,128,206,\dots$ & Unknown \\
			6   & 5 & - & $1,3,6,8,\dots$ & $50,84,150,204,\dots$ & Unknown\\
			8 & - & 5 & - & - & Non-existent \\
			10 & - & - & - & - &Non-existent \\
			12 & 7 & - & $5,7,12,14,\dots$ & $294,368,588,690,\dots$ & Unknown \\
			14 & - & 7 & - & -& Non-existent\\
			16 & - & - & - & - &Non-existent \\
			18 & - & - & - & - & Non-existent \\
			20 & 9 & - & $2,4,11,13, \dots$ & $486,580,972,1102,\dots$ & Unknown \\
			22 & - & 9 & - & - & Non-existent \\
			\hline
		\end{tabular}
	\end{center}
	\caption{The first even values for the undirected degree $r>2$ and their corresponding values for parameters $c_1$, $c_2$,  and $z$ according to Theorem \ref{th:neccond1}.}\label{tab:mixMoore}
\end{table}

In the case of diameter 2, the map $\sigma$, which assigns to each vertex $v\in V(G)$ its repeat $\sigma(v)$, is an automorphism of $G$ as in the case of digraphs of any diameter (see
Baskoro, Miller, and Plesn\'{\i}k \cite{bmp98}). From the matrix approach, the automorphism $\sigma$ can be represented as a permutation matrix $\PP$ ($p_{ij}=1$ if and only if $\sigma(i)=j$ assuming $V(G)=\{1,\ldots,n\}$). As a consequence, the adjacency matrix $\A$ of an $(r,z,2)$-almost Moore mixed graph $G$ satisfies the following matrix equation in terms of the matrix $\PP$:
\begin{equation}
	\I+\A+\A^2=\J+r\I+\PP,
	\label{eq:eqmat1}
\end{equation}
where $\I$ and $\J$ denote the identity and the all-ones matrix, respectively. A necessary condition for the existence of such graphs can be derived from Eq. \eqref{eq:eqmat1} using spectral graph theory  (see L\'opez and Miret \cite{lm16}).

\begin{theorem}[\cite{lm16}] \label{th:neccond1}
	Let $G$ be a (totally regular) almost Moore mixed graph of diameter two, undirected (even) degree $r>2$, and directed degree $z \geq 1$. Then, one of the following conditions must hold:
	\begin{itemize}
		\item[$(a)$] There exists an odd integer $c_1 \in \mathbb{Z}$ such that $c_{1}^2=4r+1$ and $c_1 \ | \ (4z+1)(4z-7)$.
		\item[$(b)$] There exists an odd integer $c_2 \in \mathbb{Z}$ such that $c_{2}^2=4r-7$ and $c_2 \ | \ (16z^2+40z-23)$.
	\end{itemize}
\end{theorem}

For instance, the possible existence of almost Moore mixed graphs of diameter $k=2$ is as in Table \ref{tab:mixMoore} for $r\in\{4,6,\ldots, 22\}$. 

\begin{figure}[t]
\begin{center}
	\includegraphics[width=14cm]{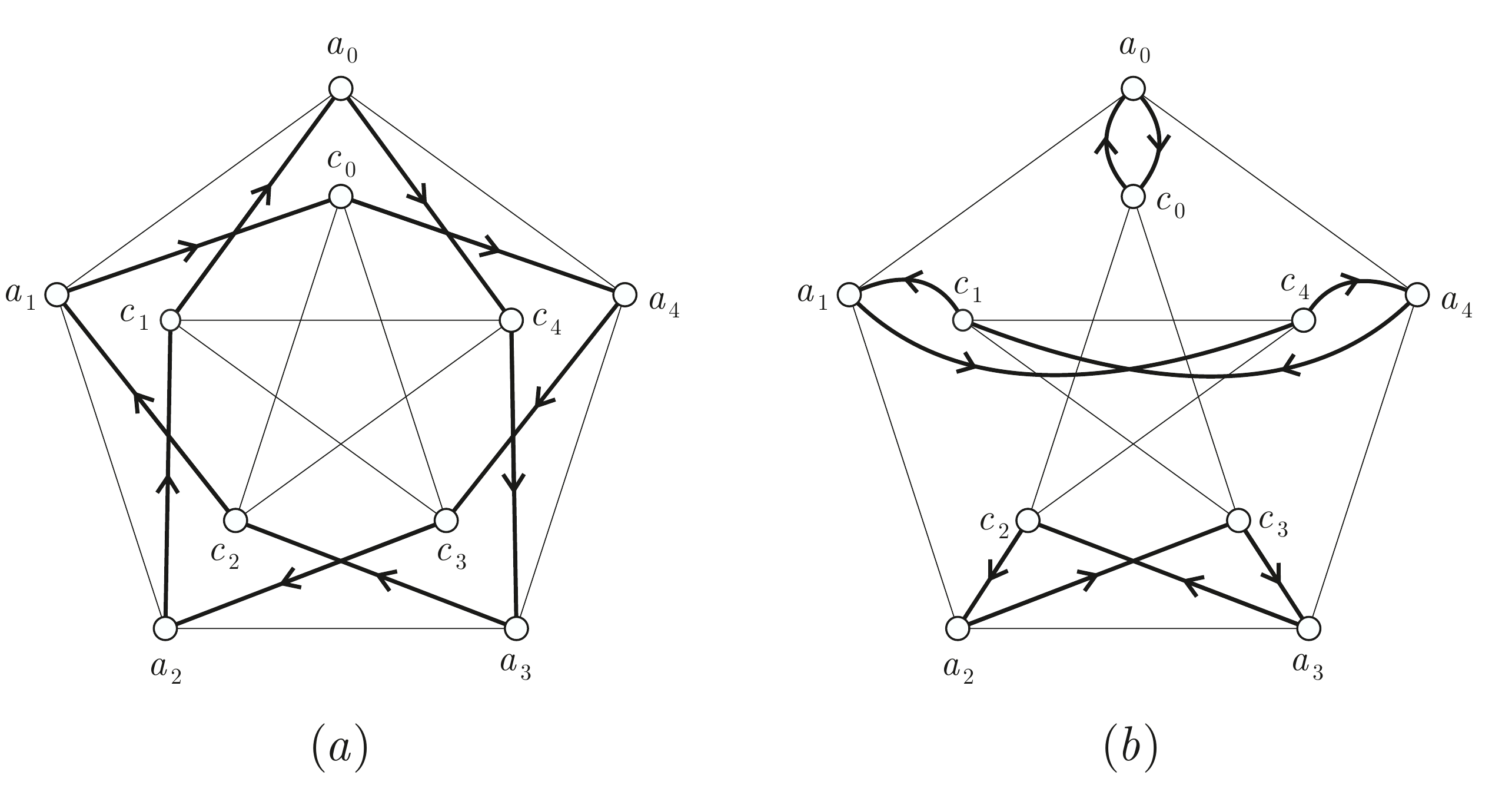}
\end{center}
\vskip-.5cm
\caption{$(a)$ The only almost Moore mixed graph of diameter $2$ known until now. $(b)$ A mixed graph of order $10$ and diameter $2$ satisfying Eq. \eqref{eq:eqmat1} that is not totally regular.}\label{almost-moore-mixed}
\end{figure}

Besides these conditions, only one almost Moore mixed graph of diameter $2$ is known until now: the one with parameters $(r,z,k)=(2,1,2)$ (see Figure \ref{almost-moore-mixed} (a)). In this particular graph, the repeat of vertex $a_i$ is $a_{i-2}$ since there are two different paths of length $\leq 2$ joining them ($a_ia_{i-1}a_{i-2}$ and $a_ic_{i-1}a_{i-2}$). Besides, $\sigma(c_i)=c_{i-2}$ since $c_ic_{i-2}$ and $c_ia_{i-1}c_{i-2}$ are again two different paths of length $\leq 2$, where all operations are considered modulo $5$. Hence, the permutation $\sigma$ decomposes in two disjoint cycles of length five, that is, $\sigma=(a_0 a_3 a_1 a_4 a_2)(c_0 c_3 c_1 c_4 c_2)$. So, the permutation cycle structure $(m_1,\dots,m_{10})$ of this graph is $m_5=2$ and $m_i=0$, for all $i \neq 5$. For more information on its existence, see L\'opez and Miret \cite{lm16}; and on its unicity, see  Buset, L\'opez, and  Miret \cite{blm17}).
Moreover, answering an open question of the two last authors \cite{lm16}, Tuite and Erskine \cite{tg19} proved that  $(r,z,2)$- and $(1,1,k)$-almost Moore mixed graphs are totally regular (that is, all vertices have undirected degree $r$ and out- and in-degree $z$).

We would like to point out that there are many mixed graphs whose adjacency matrix satisfy equation \eqref{eq:eqmat1}, but they are not totally regular in the following sense: A digon must always be considered as an edge. Otherwise, the mixed graph depicted in Figure \ref{almost-moore-mixed}(b) would be a $(2,1,2)$-almost Moore mixed  graph since vertices $a_0$ and $c_0$ would have the right degrees (instead of undirected degree $3$ and directed degree $0$, the value that they actually have). From the matrix point of view, the permutation $\sigma$ of this graph decomposes in $\sigma=(a_0)(c_0)(a_1 c_4 a_4 c_1)(a_2 c_3 a_3 c_2)$.

\section{Almost Moore mixed graphs of diameter $3$}
\label{sec:k=3}
\begin{figure}[t]
    \begin{center}
        \includegraphics[width=15cm]{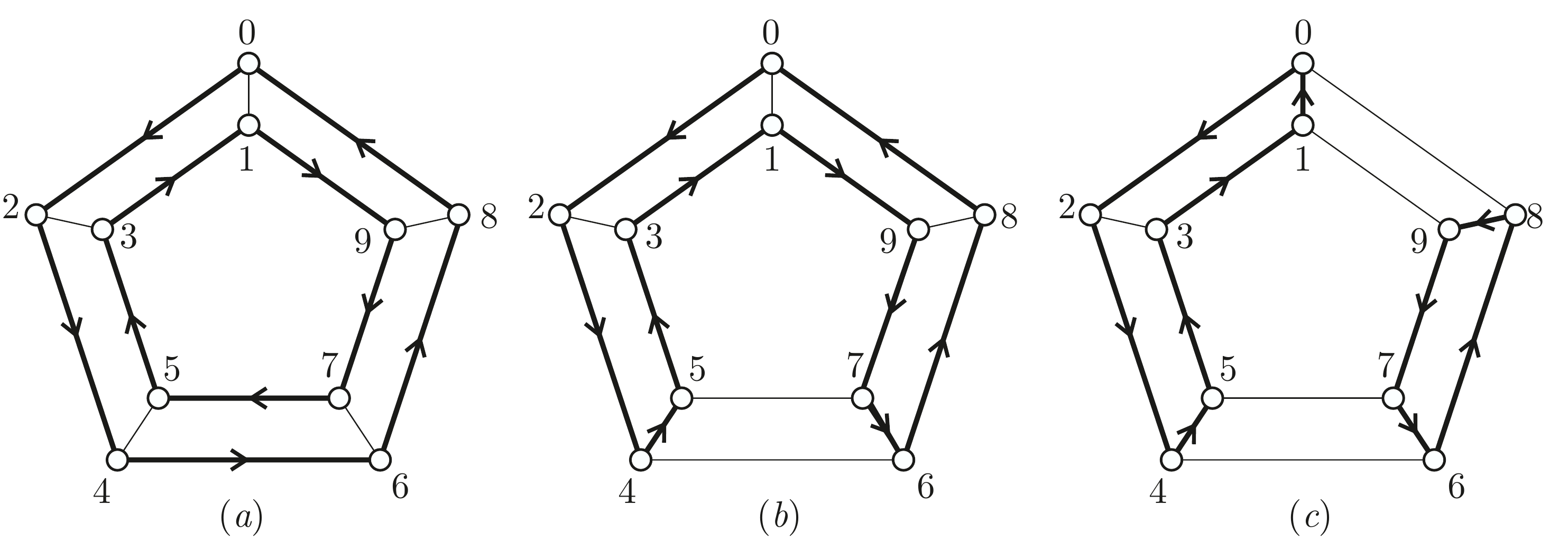}
    \end{center}
    \vskip-.5cm
	\caption{The unique three non-isomorphic almost Moore mixed graphs with diameter $k=3$ and directed degree $z=1$.}
	\label{fig1}
\end{figure}

Since Moore mixed  graphs do not exist for diameter $k \geq 3$, one could ask for the existence of almost Moore mixed graphs. In \cite{dfl18}, Dalf\'o, Fiol, and L\'opez proved the following result.
\begin{theorem}[\cite{dfl18}]
\label{theo1}
The order $N$ of an $(r,z)$-regular mixed graph $G$ with diameter $k\geq 3$ satisfies the bound
\begin{equation}
\label{improved bound}
N\le M(r,z,k)-r,
\end{equation}
where $M(r,z,k)$ is given by \eqref{eq:moorebound4}.
\end{theorem}
This means that almost Moore mixed graphs of diameter $k \geq 3$ may only exist for $r=1$ and an even $N$ (because the edges constitute a perfect matching). In fact, they showed that there exist exactly three almost Moore mixed graphs in the case $k=3$ and $z=1$ (see Figure \ref{fig1}). The graph $(a)$ depicted in Figure \ref{fig1}  has $m_2=5$ and $m_i=0$ for any $i \neq 2$, since $\sigma = (0 \ 1)(2 \ 3)(4 \ 5)(6 \ 7)(8 \ 9)$; note that, in this case, the cycles correspond to the edges. Besides, the graph $(b)$ has $m_2=3$ and $m_4=1$ (the remaining values of $m_i$ are zero) since $\sigma = ( 0 \ 1)(2 \ 3)(4\ 6 \ 7 \ 5)(8 \ 9)$. Finally, the graph $(c)$ has $m_2=1$ and $m_4=2$ since $\sigma = (2 \ 3)(4 \ 6 \ 7 \ 5)(8 \ 0 \ 1 \ 9 )$.

\begin{theorem}
\label{th:k=3}
Let $\A$ be the adjacency matrix of an almost Moore mixed graph $G$ of diameter $3$. Let $\A=\R+\Z$,
 where $\R$ and $\Z$ are  the adjacency matrices of the subgraphs of $G$ induced by the edges and arcs of $G$, respectively.
\begin{itemize}
\item[$(a)$]
There exists a permutation matrix $\PP$ such that $\A$ satisfies 
\begin{equation}
\label{k=3-Nacho}
\A^2 + \A^3 = \J + \Z + \PP.
\end{equation}
\item[$(b)$]
The mapping $\sigma$ represented by $\PP$ is not an automorphism of $G$.
\item[$(c)$]
If $\PP=\R$, then
$G$ is the Cayley graph of the dihedral group  with generators $r$ and $s$, and presentation $D_5=\langle r,s\, |\, r^5\!=\!s^2\!=\!(rs)^2\!=\!1 \rangle$.
\end{itemize}
\end{theorem}

\begin{proof}
$(a)$
Let $G$ be an almost $(r,z)$-Moore mixed  graph of diameter $3$ and order $M(r,z,3)-1$. According to Theorem \ref{theo1}, $r$ must be $1$, and so $\R$ is a permutation matrix. Now, we show that the adjacency matrix $\A$ of $G$ satisfies the matrix equation $\I+\A+\A^2+\A^3=\J+\I+\A+\Z+\PP$. 
\begin{itemize}
	\item[$(i)$] Indeed, each entry in $\I+\A+\A^2+\A^3$ is at least $1$ due to the uniqueness of the shortest walks in $G$ (up to distance $3$). This implies the $\J$ on the right side of the matrix equation.  
	\item[$(ii)$]  The fact that each vertex is incident to one edge gives exactly one closed walk (performed by edges) of length $2$ at any vertex. This corresponds to $\I$ on the right side of the equation. 
	\item[$(iii)$]  Let $u$ be the unique vertex adjacent from a fixed vertex $v$ by an edge. Then, $u$ is reached from $v$ also through the walk $v-u-v-u$ of length $3$ (corresponding to the matrix $\R$). 
	\item[$(iv)$] Let $\{u_1,u_2,\dots,u_z\}$ be the set of vertices adjacent from $v$ by an arc.  Every $u_i$ is reached from $v$ through the walks $v-u-v-u_i$ of length $3$ 
	(corresponding to the matrix $\Z$).
	 \item[$(v)$] Let $u'_i$ be the (unique) vertex adjacent from $u_i$ by an edge, for all $1 \leq i \leq z$. Every $u_i$ is reached from $v$ through the walks $v-u_i-u'_i-u_i$ of length $3$,  
	 (corresponding again to the matrix $\Z$).
	 \item[$(vi)$] Finally, every extra walk of length at most $3$ from $v$ to $\sigma(v)$ is codified in the matrix $\PP$.
\end{itemize}
Note that from $(iii)$, $(iv)$, and $(v)$, $\A^3$ counts (besides the unique shortest walks to vertices of distance $3$) one extra walk to vertices at distance $1$ (those counted in $\R+\Z=\A$) plus one extra walk to vertices at distance $1$ pointed by an arc (those counted in $\Z$). 

Altogether, we get the equation $\I+\A+\A^2+\A^3=\J+\I+\R+2\Z+\PP$. That is,
\begin{equation}\label{eq:k=3}
\A^2 + \A^3 = \J + \Z + \PP.
\end{equation}
	Under the above hypothesis on $\R$, an alternative way to obtain \eqref{eq:k=3} is to use that  $\R^{\ell}=\I$ if $\ell$ is even, and
	$\R^{\ell}=\R$ if $\ell$ is odd. Let $\A_i$ be the $i$-distance matrix of $G$, where $(\A_i)_{uv}=1$ if the distance from vertex $u$ to $v$ is $i$, and $(\A_i)_{uv}=0$ otherwise. Then, from $\A_0=\I$, $\A_1=\A$,
	\begin{align*}
		\A^2 &=(\R+\Z)^2=\I+\R\Z+\Z\R+\Z^2=\I+\A_2 \Rightarrow\\
		\A_2& =\A^2-\I, \mbox{\ and}\\
		\A^3 & =(\R+\Z)^3=\R+\Z+\R\Z\R+\R\Z^2+\Z+\Z\R\Z+\Z^2\R+\Z^3 \\
		& =\A_3+\R+2\Z =\A_3+\A+\Z \Rightarrow\\
       \A_3 &  =\A^3-\A-\Z,
	\end{align*}
	we get $\sum_{i=0}^3\A_i=\A^3+\A^2-\Z=\J+\PP$, as claimed.

$(b)$  Moreover, $\PP$ is an automorphism of $G$  if and only if $\PP$ commutes with $\A$. Then, in this case, $\Z=\A^3+\A^2-\J-\PP$ commutes with $\A$ since the regularity of $G$ implies that $\J$ is a polynomial in $\A$ (see Hoffman and McAndrew \cite{HoMc65}). Hence, $\Z$ also commutes with $\R=\A-\Z$. But the equality $\Z\R=\R\Z$ would imply that between two vertices, say $u$ and $v$, there is more than one shortest path of length two. Thus, vertex $u$ would have more than one repeat, against the hypothesis.

$(c)$ If $\PP=\R$ (that is,  the permutation $\sigma$ is involutive, and
every vertex $v$ and its repeat $\sigma(v)$ are joined by an edge), we have that
$\PP+\Z=\A$, and 
Eq. \eqref{eq:k=3} becomes 
\begin{equation}\label{eq:k=3-simpl}
-\A + \A^2+ \A^3 = \J.
\end{equation}
Note that our assumption that $\PP=\R$ is consistent with the fact that, according to Eq. \eqref{eq:k=3}, $\Z+\PP$ commutes with $\A$.

Since the spectrum of $\J$ has eigenvalues $n(=M(1,z,3)-1=(1+z)^3+(1+z)^2-(1+z))$ with multiplicity 1, and 0 with multiplicity $n-1$, $\A$ has eigenvalues $1+z$ with multiplicity 1 and the zeros of the polynomial $-x+x^2+x^3$, that is, $0$, $\alpha=\frac{\sqrt{5}-1}{2}$ and $\overline{\alpha}=\frac{-\sqrt{5}-1}{2}$.
As a consequence, the characteristic polynomial of $G$ is
\[
\phi_G(x)=(x-(1+z))x^a(x-\alpha)^b(x-\overline{\alpha})^c,
\]
where the eigenvalue multiplicites $a$, $b$, and $c$ are positive integers such that
$\tr{\A}^0=1+a+b+c=n$.
Moreover, $\tr{\A}^1=1+z+b\alpha+c\overline{\alpha}=0$ ($G$ has no loops) implies that $2(1+z)=(b+c)+(b-c)\sqrt{5}$.
Since $b,c,z$ must be positive integers, we have that $b-c=0$ and, therefore, $b=c=z+1$.
Finally, since there is one closed walk of length $2$ for any vertex in $G$, and $G$ does not contain selfrepeat vertices, we have that $\tr{\A^2}=(1+z)^2+b\alpha^2+c\overline{\alpha}^2=n$. Altogether, we get the equalities
$b=\frac{1}{3}(z(z+1)(z+2)) = z+1$, with the only feasible solution is $z=1$. Hence, $G$ has degree parameters $r=z=1$. In \cite{dfl18}, the authors proved that, in this case,  there exist only three almost Moore mixed graphs (namely, the ones depicted in Figure \ref{fig1}). Just the graph of the case $(a)$ satisfies $\sigma^2=id$, which is precisely the Cayley graph of the dihedral group $D_5=\langle r,s\, |\, r^5\!=\!s^2\!=\!(rs)^2\!=\!1 \rangle$, with generators $r$ and $s$.
\end{proof}

Notice that, as a consequence of Theorem \ref{th:k=3}$(b)$, $\PP\neq \I$, so we get the following consequence.
\begin{corollary}
	There is no $(1,z,3)$-almost Moore mixed graph $G$ with every vertex
	a selfrepeat.
\end{corollary} 

\subsection{The structure of the $(1,1,3)$-almost Moore mixed graphs}

\begin{figure}[t]
    \begin{center}
        \includegraphics[width=12cm]{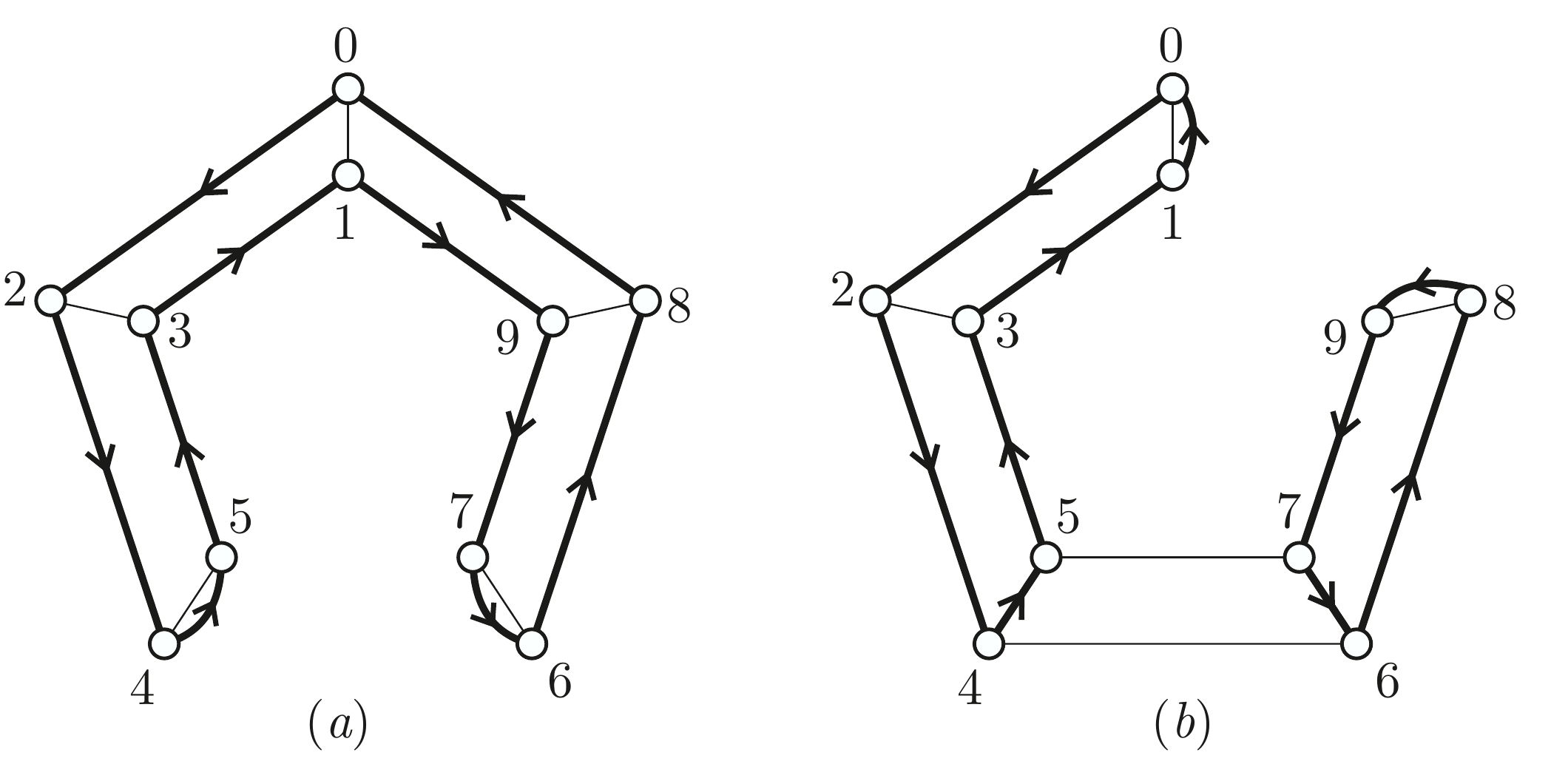}
    \end{center}
    \vskip-.5cm
	\caption{The two non-isomorphic mixed graphs with parallel arcs and cospectral with $H^{(i)}$, for $i=1,2,3$: $(a)$ $H^{(4)}(\cong H^{(6)}\cong H^{(7)})$, $(b)$ $H^{(5)}$.}
	 \label{fig2}
\end{figure}

Now let us take a close look at the structure of the only three $(1,1,3)$-almost Moore mixed graphs (see again Figure \ref{fig2}).
First, let us recall the known properties of these mixed graphs $(a)$, $(b)$, and $(c)$, from now on called $H^{(1)}$, $H^{(2)}$, and $H^{(3)}$, respectively (see Dalfó, Fiol, and López \cite{dfl18}).
Apart from being a Cayley digraph, $H^{(1)}$ is the line digraph of the cycle $C_5$ (seen as a digraph, so that each edge corresponds to a digon, that is, two opposite arcs).
The mixed graphs $H^{(2)}$ and $H^{(3)}$ can be obtained from $H^{(1)}$ by applying a recent method to obtain cospectral digraphs from a locally line digraph (see Dalfó and Fiol \cite{DaFi16}). These mixed graphs can also be obtained as a proper orientation of the pentagonal prism graph or the so-called Yutsis graph of the 15j symbol of the second kind (see Yutsis, Levinson, and Vanagas [7]). Finally, each of the three mixed graphs is isomorphic to its converse (where the directions of the arcs are reversed), and they are cospectral (see Lemma \ref{le:RZP}). 

To describe the new properties and according to the notation of the last theorem and in Figure \ref{fig2}, let $\A^{(i)}=\R^{(i)}+\Z^{(i)}$ be the adjacency matrix of $H^{(i)}$ for $i=1,2,3$. As mentioned above, let $\PP^{(i)}$ be the corresponding permutation matrices representing the permutation $\sigma^{(i)}$ (of the repeats of $H^{(i)}$). Since $r=z=1$, the matrices $\R^{(i)}$ and $\Z^{(i)}$ also are permutation matrices corresponding to the permutations, say, $\rho^{(i)}$ and $\omega^{(i)}$ respectively. More precisely,
$$
\begin{array}{lll}
	\sigma^{(1)}=(01)(23)(45)(67)(89), & \rho^{(1)}=\sigma^{(1)}, & \omega^{(1)}=(02468)(19753),\\
	\sigma^{(2)}=(01)(23)(4675)(89), & \rho^{(2)}=(01)(23)(57)(46)(89), & \omega^{(2)}=(0245319768),\\
	\sigma^{(3)}=(23)(4675)(8019), & \rho^{(3)}=(08)(23)(57)(46)(19), & \omega^{(3)}=(024531)(6897).
\end{array}
$$
From this, it is routine to check the following lemma describing the new properties.
\begin{lemma}
\label{le:RZP}
Let $\PP^{(i)}$, $\R^{(i)}$, and $\Z^{(i)}$ the matrices corresponding to the above permutations $\sigma^{(i)}$, $\rho^{(i)}$, and  $\omega^{(i)}$, respectively. Let $\A^{(4)}=\R^{(1)}+\Z^{(2)}$, $\A^{(5)}=\R^{(2)}+\Z^{(3)}$, $\A^{(6)}=\PP^{(1)}+\Z^{(2)}$, and $\A^{(7)}=\PP^{(2)}+\Z^{(3)}$, which correspond to adjacency matrices of the mixed graphs with parallel arcs $H^{(4)},H^{(5)},H^{(6)},H^{(7)}$, respectively (see Figure \ref{fig2}). Then,
\begin{itemize}
	\item[$(a)$]
 $(\R^{(i)})^{2}=\I$, and $\A^{(i)}\A^{(j)}=\A^{(j)}\A^{(i)}$ for all $i,j\in \{1,2,3\}$.
	\item[$(b)$]
$\A^{(1)}=\PP^{(i)}+\Z^{(i)}$ for all $i\in \{1,2,3\}$.
	\item[$(c)$]
$\A^{(i)}=\R^{(i)}+\Z^{(i)}=(\PP^{(i)})^{\top}+\Z^{(1)}$ for all $i\in \{1,2,3\}$.
	\item[$(d)$]
All the mixed graphs $H^{(i)}$ are cospectral, with spectrum
$$\textstyle \spec (\A^{(i)})=\left\lbrace 2^{(1)}, 0^{(5)}, \left(\frac{-1\pm\sqrt{5}}{2}\right)^{(2)}\right\rbrace $$
for  $i=1,\ldots,7$. (In fact, $H^{(4)}$, $H^{(6)}$ and $H^{(7)}$ are isomorphic.)
\end{itemize}
\end{lemma}

Note that, from $(b)$, $\PP^{(i)}=\R^{(1)}+\Z^{(1)}-\Z^{(i)}$. Moreover, from $(c)$, $\R^{(i)}=(\R^{(1)})^{\top}+(\Z^{(1)})^{\top}-(\Z^{(i)})^{\top}+\Z^{(1)}-\Z^{(i)}$. Then, using $(a)$ and reordering the terms, we get
$$
\R^{(i)}+[\Z^{(i)}+(\Z^{(i)})^{\top}]=\R^{(1)}+[\Z^{(1)}+(\Z^{(1)})^{\top}] \quad \mbox{ for all } i\in \{1,2,3\}.
$$
Both sides of this equation correspond to the adjacency matrix of the underlying graph of $H^{(i)}$, that is, the pentagonal prism graph for $i\in \{1,2,3\}$, as commented in the beginning of this subsection. 

As another 
consequence of this lemma, we have the following result.

\begin{proposition}
Let us consider what we called the {\em $(1,1,3)$-mixed adjacency algebra} of matrices ${\cal A}_{(1,1,3)}=\langle \R^{(1)},\Z^{(1)},\Z^{(2)},\Z^{(3)} \rangle$. Then,
$$
\A^{(i)}\in {\cal A}_{(1,1,3)},\qquad i=1,\ldots, 7.
$$
\end{proposition}
\begin{proof}
More generally, we prove that $\R^{(i)},\Z^{(i)},\PP^{(i)}\in {\cal A}_{(1,1,3)}$.
First, by using that $\A^{(1)}=\R^{(1)}+\Z^{(1)}$  and Lemma \ref{le:RZP}$(b)$,
we get that $\PP^{(i)}=\A^{(1)}-\Z^{(i)}=\R^{(1)}+\Z^{(1)}-\Z^{(i)}$ for $i=1,2,3$ (in particular, $\PP^{(1)}=\R^{(1)}$ and, hence, $(\PP^{(1)})^{2}=\I$).
Now, from Lemma \ref{le:RZP}$(c)$, $\R^{(i)}=(\PP^{(i)})^{\top}+\Z^{(1)}-\Z^{(i)}$ for $i=1,2,3$ (the case $i=1$ is trivial since $(\PP^{(1)})^{\top}=\PP^{(1)}$).
\end{proof}

\section{Large mixed graphs from line digraphs}
\label{sec:construc}
Some constructions of large mixed graphs have been proposed in the literature. For example, the best infinite families of such graphs with an asymptotically optimal number of vertices for their diameter were proposed by Dalf\'o in  \cite{d19}. These graphs are vertex-transitive and
generalize both the pancake graphs \cite{Dw} (when $z= 0$) and the Faber-Moore-Chen or cycle-prefix digraphs \cite{FaMoCh} (when $r= 1$).

Here we propose a construction based on the line digraph technique.
Given a digraph $G=(V,E)$, its  line digraph $LG$ has vertices representing
the arcs of $G$, and each vertex $uv$ (where $u\rightarrow v$ in $G$) is adjacent to the vertices $vw$ for all $w$ adjacent from $v$ in $G$. As already  commented, Figure \ref{fig1}$(a)$ 
shows the line digraphs of the cycle $C_5$,
with edges corresponding to digons.
Our first construction is based on the following result by Fiol, Yebra, and Alegre   \cite{FiYeAl84}, where the average distance of $G$ is $\overline{k}=\frac{1}{n^2}\sum_{u,v\in V} \dist(u,v)$.

\begin{theorem}[\cite{FiYeAl84}]
Let $G$ be a $\delta$-regular digraph $(\delta>1)$ of order $n$, diameter $k$, and average distance $\overline{k}$. Then, the order $n_L$, diameter $k_L$, and average distance $\overline{k}_L$ of the line digraph $LG$ satisfy
\begin{equation}
\label{n,k-LH}
n_L=\delta n,\quad k_L=k+1,\quad\mbox{and}\quad \overline{k}_L<\overline{k}+1.
\end{equation}
\end{theorem}

\begin{lemma}
For every $\delta$-regular graph $H$  with $n$ vertices, diameter $k$, and average distance $\overline{k}$, there is a $(1,\delta-1,k+1)$-mixed graph $G$ with undirected degree $r=1$, directed degree $z=\delta-1$, order $\delta n$, diameter $k+1$, and average distance smaller than $\overline{k}+1$.
\end{lemma}

\begin{proof}
Consider the digraph $H'$ obtained from the graph $H$, where each digon of $H'$ corresponds to an edge of $H$. Then, the line digraph $G=LH'$ has the claimed parameters  since each edge of $H$ (that is, each digon $u\rightarrow v$ and $v\rightarrow u$ of $H'$) gives rise to one edge (digon $uv\rightarrow vu$ and $vu\rightarrow uv$) of $G$, so $r=1$. Moreover, each arc $u\rightarrow v$ of $H'$ is adjacent to $\delta-1$ arcs $v\rightarrow w$, with $w \neq u$. Hence, vertex $uv$ of $G$ is adjacent to $\delta-1$ vertices $vw$ of $G$, so $z=\delta-1$. The order, diameter, and average distance of $G$ follow from \eqref{n,k-LH}.
\end{proof}
Some examples of  large mixed graphs obtained by applying this lemma follow:
\begin{itemize}
\item
The mixed graph $LK_n$, with $n=d+1$,  is an almost Moore mixed graph with $r=1$, $z=d-1$, $N=d^2+d$ vertices, and diameter $k=2$. This is isomorphic to the well-known  Kautz digraph $K(d,2)$ and, as it was proved by Gimbert \cite{g01}, it is also the unique almost Moore digraph of diameter two.
\item
The mixed graph $LC_n$, with $n=2\ell+1$, is a mixed graph with $r=z=1$, $N=2n$ vertices, and diameter $k=\ell$. Note that, in particular, $LC_3\cong LK_3$. Moreover, as we already showed in Section \ref{sec:k=3},
$LC_5$ (see again Figure  \ref{fig1}$(a)$) is a $(1,1,3)$-almost Moore mixed graph.
\item
Let $G_1$, $G_2$, and $G_3$ be the known Moore graphs with diameter two. Namely, $G_1=C_5$, $G_2=P$ (the Petersen graph), and $G_3=HS$ the Hoffman-Singleton graph. Then, $LG_1\cong LC_5$; $LG_2$ is a $(1,2,3)$-mixed graph with $N=30$ vertices (the Moore bound is $M(1,2,3)=34$, but we know that the maximum is $32$ (the even number smaller than $M(1,2,3)-1$);
$LG_3$ is a $(1,6,3)$-mixed graph with $N=350$ vertices (the Moore bound is $M(1,6,3)=386$, but the maximum must be $384$).
We conjecture that all these mixed graphs have the maximum possible order.
\item
The mixed graph $LK_{n,n}$ is a bipartite mixed graph with $r=1$, $z=n-1$, $N=2n^2$ vertices, and diameter $k=3$. This is a Moore bipartite mixed graph because, in the case of bipartite $(1,z,3)$-mixed graphs, the corresponding Moore bound 
(given in Dobravec and Robi\v{c} \cite{dr09}) 
$$
M_b(1,z,3)=2(1+2z+z^2)=2(z+1)^2
$$
is attained. 
\end{itemize}



\end{document}